\documentclass[12pt,psfig,reqno]{amsart}
\usepackage{amssymb,amsfonts,latexsym}
\usepackage{graphics,verbatim}
\usepackage{graphicx}
\usepackage[usenames]{color} %\color{Red}
\usepackage{soul} %\textst{word}

\hoffset=-1cm \errorcontextlines=0 \numberwithin{equation}{section}

 \theoremstyle{plain}
\newtheorem{theorem}{Theorem}[section]
\newtheorem{lemma}[theorem]{\bf{Lemma}}
\newtheorem{prop}[theorem]{\bf{Proposition}}

\newenvironment{pf}{{\noindent \bf Proof.\/}}{\hfill$\Box$}

\date {}

\begin{document}

\title{Kato's inequality and Liouville theorems on locally finite graphs}

\author{Li Ma, Xiangyang Wang}

\address{Li Ma, Department of mathematics \\
Henan Normal university \\
Xinxiang, 453007 \\
China}

\email{nuslma@gmail.com}

\address{X.Y.Wang, Department of mathematics \\
Sun Yah Sen university \\
Guangzhou,\\
China}

\email{mcswxy@mail.sysu.edu.cn}

\thanks{ The research is partially supported by the National Natural Science
Foundation of China 10631020 and SRFDP 20090002110019}

\begin{abstract}
In this paper we study the Kato' inequality on locally finite graph.
We also study the application of Kato inequality to Ginzburg-Landau
equations on such graphs. Interesting properties of Schrodinger
equation and a Liouville type theorem are also derived.

{ \textbf{Mathematics Subject Classification 2000}: 31C20, 31C05}

{ \textbf{Keywords}: locally finite graph, Kato's inequality,
Ginzburg-Landau equation, Liouville theorem}
\end{abstract}

\maketitle

\section{Introduction}
In recent studies Yau and F.Chung and their friends (see \cite{Ch},
\cite{CY} and \cite{LY} for more background and references) have
studied Ricci curvature and eignevalue estimate on locally finite
graphs. The lower bound of Ricci curvature on locally finite graphs
can be defined via the method of Bakry-Emery. Then following the
method of Li-Yau, one can do the gradient estimate for
eigen-functions of the Laplacian operators on locally finite graphs.
In particular, one can derive the lower bound of the eigenvalues on
a connected graph with finite diameter. On a connected graph with
finite diameter, one can see that the Liouville theorem is always
true for harmonic function. In fact, the harmonic functions are
always bounded in the connected finite graphs and their maximum
values are obtained somewhere. Using the mean value property, one
obtain the Liouville theorem. It is nature to ask if such a
Liouville type theorem is true on nonlinear elliptic problems on
locally finite graphs. With this question in mind, we intend to
study the Kato's inequalities in this paper. As an application we
get a Liouville theorem for nonlinear elliptic equations on the
locally finite graphs. Our results are stated in lemmas
\ref{lemma-Kato-1} and \ref{lemma-kato-2} below.

We mention the other motivation of this paper. The Ginzbourg-Landau
equation is a basic model for the mathematical theory of
superconductivity, which examines the macroscopic properties of a
superconductor with the aid of general thermodynamic arguments
\cite{LP}. This equation is derived from the free energy of the form
$$
\frac{1}{2}|\nabla u|^2+\frac{1}{4}(1-|u|^2)^2
$$
of the complex order parameter $u$. We shall confine the complex
variable $u$ defined on locally finite graphs $X$ and study the
property of the solutions of the Ginzburg-Landau equation
$$
-\Delta u+u(|u|^2-1)=0, \ \ \ in \ \ X.
$$
With the help of Kato's inequality we show the uniform bound of the
solutions $u$ such that $|u|\leq 1$ on $X$. Related works on the
whole Euclidean space can be found in \cite{M} and \cite{MX}.

We also study the interesting properties of Schrodinger equations on
locally finite graphs.

The plan of the paper is below. Notations are introduced in section
\ref{sect2} and all of our results are stated and proved in section
\ref{sect2}.

\section{Set up and proofs of main results}\label{sect2}

Let $(X, {\mathcal E})$ be a graph with countable vertice set $X$
and edge set ${\mathcal E}$.  We assume that the graph is {\it
simple}, i.e., no loop and no multi-edges. We also assume that the
graph is connected. Let $\mu_{xy}= \mu_{yx} >0$ is a symmetric
weight on ${\mathcal E}$. We call $d_x = \sum_{(x,y)\in{\mathcal E}}
\mu_{xy}$
(we also assume $d_x<\infty$ for all $x\in X$) the {\it degree of } $x\in X$.\\

Denote by
\[
\ell(X) = \{ u:\ u:\ X \longrightarrow {\Bbb R}\},
\]
the set of all real functions (or complex-valued functions with
$\Bbb R$ replaced by $\Bbb C$ on $X$. We often denote by $u^2$ as
$|u|^2$.

We define the {\it Laplacian operator} $\Delta: \ \ell(X) \longrightarrow \ell(X)$:
\[
(\Delta u) (x) = \sum_{(x,y)\in {\mathcal E}} \frac{\mu_{xy}}{d_x}\big( u(y) - u(x)\big).
\]
We also define
\[
|\triangledown u|^2 (x) = \sum_{(x,y)\in {\mathcal E}} \frac{\mu_{xy}}{d_x}
\big(u(y) - u(x)\big)^2.
\]

Then we have the following elementary fact.

\begin{lemma} \label{lemma-Kato-1}
(Kato's inequality) For a graph $X$, we have
\[
| \triangledown u|^2 \ge \big|\triangledown |u| \big|^2.
\]
\end{lemma}

\begin{pf}
For any $x\in X$, we have
\[
\begin{aligned}
|\triangledown u|^2 (x) & = \sum_{(x,y)\in {\mathcal E}} \frac{\mu_{xy}}{d_x} \big( u(y) - u(x)\big)^2 \\
& \ge \sum_{(x,y)\in{\mathcal E}} \frac{\mu_{xy}}{d_x} \big( | u(y)| - | u(x) |\big)^2 = |\triangledown \big| u |\big|^2 (x).
\end{aligned}
\]
This completes the proof.
\end{pf}

Generally speaking, given $\Delta u$, one may not have the
well-defined $\Delta u^2$ on fractals. However, this is not the case
on graphs.
\begin{lemma}\label{lemma-2}
$\Delta u^2 = 2 u \Delta u + |\triangledown u|^2.$
\end{lemma}

\begin{pf}
For any $x\in X$,
\[
\begin{aligned}
(\Delta u^2) (x) & = \sum_{(x,y)\in {\mathcal E}} \frac{\mu_{xy}}{d_x} \big( u^2(y) - u^2(x) \big)\\
& = \sum_{(x,y)\in {\mathcal E}} \frac{\mu_{xy}}{d_x} \big( 2 u(x) ( u(y) -u(x)) + (u(y) - u(x))^2 \big)\\
& = 2 u(x) \Delta u(x) + |\triangledown u |^2 (x).
\end{aligned}
\]
\end{pf}

With the help of above fact, we have
\begin{lemma}\label{lemma-kato-2}
(Kato's inequality)
\begin{eqnarray}
\Delta |u| & \ge {\rm sign} (u) \Delta u, \label{eq-kato-2}\\
\Delta u_+ & \ge {\rm sign}_+(u) \Delta u.\label{eq-kato-3}
\end{eqnarray}

\end{lemma}

\begin{pf} By Lemma \ref{lemma-2}, we have
\[
\Delta u^2 = 2 u \Delta u + | \triangledown u |^2
\]
and
\[
\Delta u^2 = \Delta |u|^2 = 2 |u| \Delta |u| + \big|\triangledown |u|\big|^2
\]
Hence
\[
2|u| \Delta |u| = 2 u \Delta u + |\triangledown u|^2 - \big| \triangledown |u|\big|^2.
\]
By Lemma \ref{lemma-Kato-1}, we have
\[
|u| \Delta |u| \ge u \Delta u
\]
It follows that
\[
\Delta |u| \ge \frac{u}{| u |} \Delta u = {\rm sign }(u) \Delta u,
\]
providing $u(x) \not= 0$. If $u(x)=0$, then
\[
\Delta u(x) = \sum_{(x,y)\in{\mathcal E}} \frac{\mu_{xy}}{d_x} u(y) \le \sum_{(x,y)\in{\mathcal E}} \frac{\mu_{xy}}{d_x} |u(y)| = \Delta |u| (x).
\]
we see that (\ref{eq-kato-2}) still hold.

To prove (\ref{eq-kato-3}), we note that $u_+ = \frac{1}{2} (|u| + u)$, hence
\[
\begin{aligned}
\Delta u_+ & = \frac{1}{2} ( \Delta |u| + \Delta u)\\
& \ge \frac{1}{2} \big( {\rm sign}(u) \Delta u + \Delta u \big)\\
& = \frac{1}{2} \big( {\rm sign}(u) + 1 \big) \Delta u\\
& = {\rm sign}_+(u) \Delta u.
\end{aligned}
\]
This completes the proof.
\end{pf}

We now use Kato's inequality to study properties of solutions to the
Ginzburg-Landau equation on graphs.

\begin{theorem}
Assume that $u$ is a solution of the following Ginzburg-Landau
equation
\[
\Delta u + u(1-u^2) = 0, \ \ in \ \ X.
\]
Then $|u| \le 1$.
\end{theorem}
\begin{pf}
Let $w = u^2 -1$, then
\[
\begin{aligned}
\Delta w & = 2 u \Delta u + |\triangledown u |^2\\
& = 2 u \cdot u (u^2 -1) + |\triangledown u|^2\\
& = 2(w+1)w + |\triangledown u|^2.
\end{aligned}
\]
Hence
\[
\begin{aligned}
\Delta w_+ & \ge {\rm sign}_+(w) \Delta w\\
& \ge {\rm sign}_+(w) \big( 2(w+1)w + |\triangledown u|^2 \big)\\
& \ge 2 w_+^2 + 2w_+.
\end{aligned}
\]
Assume that $\varphi > 0$ such that $ - \Delta \varphi = \lambda \varphi$ for some $\lambda>0$, then
\[
0\le \int 2(w_+^2 + 2 w_+) \varphi \le \int \varphi \Delta w_+ = \int (\Delta \varphi) w_+ = - \int \lambda \varphi w_+ \le 0.
\]
It follows that $w_+ = 0$, i.e., $w \le 0$. Hence $u^2 \le 1$.
\end{pf}

\begin{prop}
Assume $Q \ge 0\in \ell(X)$ and let $u$ be a solution such that
\begin{equation}\label{eq-Q}
-\Delta u + Q u = 0.
\end{equation}
Then $u_+$ is a sub-solution of (\ref{eq-Q}).
\end{prop}

\begin{pf}
 By the Kato's inequality, we have
 \[
 \Delta u_+ \ge {\rm sign}_+(u) \Delta u = {\rm sign}_+(u) Q u = Q u_+
 \]
 i.e., $-\Delta u_+ + Q u_+ \le 0$. That is to say, $u_+$ is a sub-solution to (\ref{eq-Q}).
 \end{pf}

with this understanding, we can do the gradient estimate for
solutions to the (stationary) Schrodinger equation and our result
extends slightly the gradient estimate in \cite{LY}.
\begin{theorem}
Assume that $u,\ Q \in \ell(X)$, $u \ge 0,\ Q \ge 0$, such that
$-\Delta u + Q u = 0$. Then
\[
|\triangledown u|^2(x) \le \left(d(1+Q(x))^2 - 2 Q(x) -1 \right) u^2(x) \le d Q^2(x) u^2(x), \quad \forall x\in X,
\]
where the constant $d=\sup_{x\in X} \sup_{(x,y)\in {\mathcal E}} \frac{d_x}{\mu_{xy}}$.
\end{theorem}

\begin{pf}
Observe that
\[
\Delta u (x) = \sum_{(x,y)\in {\mathcal E}} \frac{\mu_{xy}}{d_x} \big( u(y) - u(x) \big)
= \sum_{(x,y)\in {\mathcal E}} \frac{\mu_{xy}}{d_x} u(y) - u(x).
\]
Hence
\[
\sum_{(x,y)\in{\mathcal E}} \frac{\mu_{xy}}{d_x} u(y) = \Delta u(x) + u(x).
\]
By definition,
\[
\begin{aligned}
|\triangledown u|^2 (x) & = \sum_{(x,y)\in{\mathcal E}} \frac{\mu_{xy}}{d_x} \big( u(y) -u(x) \big)^2 \\
& = \sum_{(x,y)\in{\mathcal E}} \frac{\mu_{xy}}{d_x} \left( -2u(x) \big( u(y) - u(x) \big) -u^2(x) + u^2(y) \right)\\
& = -2 u(x) \sum_{(x,y)\in{\mathcal E}} \frac{\mu_{xy}}{d_x} \big( u(y) - u(x) \big) -u^2(x) + \sum_{(x,y)\in{\mathcal E}} \frac{d_x}{\mu_{xy}} \left(\frac{\mu_{xy}}{d_x} u(y) \right)^2\\
& \le -2 u(x) \Delta u(x) - u^2(x) + d \left(\sum_{(x,y)\in{\mathcal E}} \frac{\mu_{xy}}{d_x} u(y) \right)^2 \\
& = -(2 Q(x) + 1) u^2(x) + d\big( \Delta u(x) + u(x) \big)^2\\
& = -(2 Q(x) +1) u^2(x) + d \big( 1 + Q(x) \big)^2 u^2(x)\\
& = \left( d(1+ Q(x) )^2 - 2Q(x) -1 \right) u^2(x)\\
& \le d Q^2(x) u^2(x).
\end{aligned}
\]
 In the first inequality,
 we have uses $\frac{\mu_{xy}}{d_x} u(y) \ge 0$ for all $y\in X$ such
that $(x,y)\in{\mathcal E}$.
\end{pf}

We now derive the Liouville theorem along the line of the
Keller-Osserman theory.
\begin{theorem}
Assume that  $\ u \in \ell(X)$ and $0 \le u \le A$ (where $A$ is a
positive constant). If $\Delta u \ge u^p$ for some $p \in {\Bbb
R}_+$, then $u=0$.
\end{theorem}
\begin{pf}
Suppose otherwise, then there exists $x_0\in X$ such that $0< u(x_0)
:= \rho$. We let $w = \frac{u}{\rho}$. Then $0 \le w \le
\frac{A}{\rho}$, $w(x_0) = 1$ and
\[
\Delta w = \frac{1}{\rho} \Delta u \ge \frac{1}{\rho} u^p = \rho^{p-1} w^p.
\]
It follows that, for any $x\in X$,
\[
\sum_{(x,y)\in{\mathcal E}} \frac{\mu_{xy}}{d_x} \big( w(y) - w(x) \big) \ge \rho^{p-1} w^p(x).
\]
i.e.,
\[
\sum_{(x,y)\in{\mathcal E}} \frac{\mu_{xy}}{d_x} w(y) \ge w(x) + \rho^{p-1} w^p(x).
\]
Note that the left hand of the above is the (weighted) average of
$w(y)$'s. Hence there exists $y$ with $(x,y)\in {\mathcal E}$ such
taht
\[
w(y) \ge w(x) + \rho^{p-1} w^p(x).
\]
Using this and by induction, we get a sequence $\{x_n\}_{n=0}^\infty
\subset X$ with $(x_i, x_{i+1}) \in {\mathcal E}, \ 1=0,1,\cdots $
such that
\begin{equation}\label{eq-limit}
w(x_{n+1}) \ge w(x_n) + \rho^{p-1} w^p(x_n).
\end{equation}
It follows that $\{w(x_n)\}_n$ is a increasing sequence and bounded
by constant $\frac{A}{\rho}$. Hence there is a finite limit. Taking
the limit at the both side of (\ref{eq-limit}), we get
$\lim_{n\to\infty} w(x_n) = 0$. This contradicts that the sequence
is increasing and $w(x_0)=1$. This completes the proof.
\end{pf}

We have the following strong maximum principle for the Laplacian
equations on the locally finite graph $X$.
\begin{prop}
Assume that $u:X\to \Bbb R$ satisfies $\Delta u \ge 0$. If there
exists $x_0 \in X$ such that $u(x_0) = \sup_{x\in X} u(x) < \infty$,
then $u$ is a constant on $X$.
\end{prop}
\begin{pf} By the hypothesis on Laplacian, we have
\[
u(x_0) \le \sum_{(x,y)\in {\mathcal E}} \frac{\mu_{xy}}{d_x} u(y) \le u(x_0).
\]
Hence $u(y) = u(x_0)$ for all $y$ such that $(x,y) \in {\mathcal
E}$. By induction and the connectivity, we have $u(y) = u(x_0)$ for
all $y\in X$.
\end{pf}

Using a similar argument we have
\begin{prop}
Assume that $u:\ X\times [0,T] \longrightarrow {\Bbb R}$ such that $u_t = \Delta u$ and $u(x_0, t_0) = \sup\{ u(x,t):\ (x,t)\in X \times [0,T]\},\ t_0>0$, then $u$ is constant.
\end{prop}
\begin{pf}
At $(x_0, t_0)$, we have $u_t(x_0, t_0) \ge 0$. Similar as above argument, we have
$u(y, t_0) = u(x_0, t_0)$ for all $y\in X$ such that $(x_0,y)\in {\mathcal E}$. Also we have $u_t(y, t_0) \ge 0$. Repeating  this argument, we see the assertion holds.
\end{pf}

We shall see that the mass and energy conservation laws can also be
derived for the Schrodinger equations.
\begin{theorem} Assume that the initial data $u_0$ has finite $L^2$
norm $||u_0||_{L^2}$ and finite Dirichlet energy $||\triangledown
u_0 ||^2_{L^2}$. Then there is a unique solution $u:\ X\times
[0,+\infty):\ \longrightarrow {\Bbb C}$ to the Schrodinger equation
on the locally finite graph $X$:
\[
i u_t + \Delta u = 0; \quad u|_{t=0} =u_0.
\]
Then
\[
||u(t)||^2_{L^2} = ||u_0||^2_{L^2}, \quad ||\triangledown u(t) ||^2_{L^2} = || \triangledown u_0 ||^2_{L^2}, \quad t \ge 0.
\]
\end{theorem}

\begin{pf} We remark that the existence part of the solution to the Schrodinger equation is by now standard
and it can be derived as in the case of heat equation
via the fundamental solution. Hence we may omit the detail.

We denote by $(u,v)=u\cdot \bar{v}$ for complex valued
functions. Then we have $|u|^2=u\cdot \bar{u}$.

Compute directly and we have
\[
\begin{aligned}
\frac{d}{d t} ||u(t) ||^2_{L^2} & = 2 Re (u, u_t)\\
& = - 2 i Im (u, i u_t)\\
& = 2 i Im(u, \Delta u)\\
& = -2 i Im(\triangledown u, \triangledown u) =0.
\end{aligned}
\]
Similarly, we have
\[
\begin{aligned}
\frac{d}{d t} ||\triangledown u||^2_{L^2} & = 2 Re (\triangledown u, \triangledown u_t)\\
& = 2 Re(\Delta u, u_t)\\
& = 2 Re(\Delta u, i \Delta u) = 0.
\end{aligned}
\]

Hence the proof is complete.

\end{pf}

Similar result is true for the Gross-Pitaevskii equation on the
finite graph $X$:
$$
i u_t + \Delta u = u(|u|^2-1); \quad u|_{t=0} =u_0
$$
with the energy replaced by the free energy
$$
\frac{1}{2}|\nabla u|^2+\frac{1}{4}(1-|u|^2)^2
$$
and with finite free energy at initial time.

{\bf Acknowledgement}. This work is done while both authors visiting
the Department of Mathematics, CUHK, Hongkong and the authors would
like to thank the hospitality of the Mathematical Department of
CUHK.

{20}
\end{document}